\documentclass[12pt]{amsart}

\usepackage{amsmath,amsfonts,amscd,amsthm,amssymb,latexsym}

\font\teneufm=eufm10 \font\seveneufm=eufm7 \font\fiveeufm=eufm5
\newfam\frakturfam
\textfont\frakturfam=\teneufm \scriptfont\frakturfam=\seveneufm
\scriptscriptfont\frakturfam=\fiveeufm
\def\frak{\fam\frakturfam}

\newtheorem{pr}{Proposition}

\newtheorem{lm}{Lemma}
\newtheorem{theor}{Theorem}
\newtheorem{co}{Corollary}
\newtheorem{rem}{Remark}

\def\bee{\begin{eqnarray}}
\def\bes{\begin{eqnarray*}}
\def\eee{\end{eqnarray}}
\def\ees{\end{eqnarray*}}

\def\a{\alpha}
\def\b{\beta}

\def\s{\sigma}

\def\Proof{{\sl Proof.}\ }

\title{A Dixmier theorem for Poisson enveloping algebras}

\begin{document}
\date{}
\maketitle

\begin{center}

{\bf Ualbai Umirbaev}\footnote{
 Wayne State University,
Detroit, MI 48202, USA, and Institute of Mathematics and Modeling, Almaty, Kazakhstan,
Partially supported by the grant AP05133009 of MES RK,
e-mail: {\em umirbaev@math.wayne.edu}}
and
{\bf Viktor Zhelyabin}\footnote{Institute of Mathematics of the SB of RAS, Novosibirsk, 630090, Russia,
e-mail: {\em vicnic@math.nsc.ru}}

\end{center}

\begin{abstract} We consider a skew-symmetric $n$-ary bracket on the polynomial algebra $K[x_1,\ldots,x_n,x_{n+1}]$ ($n\geq 2$) over a field $K$ of characteristic zero
defined by $\{a_1,\ldots,a_n\}=J(a_1,\ldots,a_n,C)$, where $C$ is a fixed element of $K[x_1,\ldots,x_n,x_{n+1}]$ and $J$ is the Jacobian. If $n=2$ then this bracket is a Poisson bracket and if $n\geq 3$ then it is an  $n$-Lie-Poisson bracket on $K[x_1,\ldots,x_n,x_{n+1}]$.  We describe the center of the corresponding $n$-Lie-Poisson algebra and show that the quotient algebra $K[x_1,\ldots,x_n,x_{n+1}]/(C-\lambda)$, where $(C-\lambda)$ is the ideal generated by $C-\lambda$, $0\neq \lambda \in K$, is a simple central $n$-Lie-Poisson algebra if $C$ is a homogeneous polynomial that is not a proper power of any nonzero polynomial. This construction includes the quotients $P(\mathrm{sl}_2(K))/(C-\lambda)$ of the Poisson enveloping algebra $P(\mathrm{sl}_2(K))$ of the simple Lie  algebra $\mathrm{sl}_2(K)$, where $C$ is the standard Casimir element of $\mathrm{sl}_2(K)$ in $P(\mathrm{sl}_2(K))$. It is also proven that the quotients $P(\mathbb{M})/(C-\lambda)$ of the Poisson enveloping algebra $P(\mathbb{M})$ of the exceptional simple seven dimensional Malcev algebra $\mathbb{M}$ are central simple.
\end{abstract}


\noindent


\section{Introduction}

\hspace*{\parindent}

 Many interesting and important results have been obtained about
the structure of simple associative and Lie algebras.
Non-commutative Poisson algebras were studied by I. Herstein and his students  \cite{Herstein}. It is known if $R$ is a non-commutative prime Poisson algebra then the Poisson bracket must coincide with the commutator up to a multiple of an element from the extended centroid \cite{FL}.

 Nowadays the definition of Poisson algebras includes commutativity. Commutative Poisson algebras  are very important in many branches of mathematics and
physics and are investigated by many scientists.
The structure of simple commutative Poisson algebras has not been studied at all. As far as we know, there is only one paper  \cite{Farkas} devoted to the structure of simple Poisson algebras. The main examples of simple algebras are the symplectic Poisson algebras $P_n$ (see Section 2). A symplectic Poisson algebra by Farkas \cite{Farkas} is a Poisson algebra $P$ such that $\mathrm{Der}(P)=P \cdot\mathrm{Ad}(P)$ (see also \cite{FL}), i.e., the derivations of the commutative associative algebra $P$  as a $P$-module  are generated by all  inner Poisson derivations of $P$. Farkas proved \cite{Farkas} that commutative regular affine symplectic domains over a field of characteristic zero are Poisson simple and gave an example of a  simple non-symplectic Poisson algebra on polynomial algebras in three variables.

This paper was initiated by the following amazing result by Dixmier  \cite{Dixmier73}. In 1973 he studied the quotients $U(\mathrm{sl}_2(K))/(C-\lambda)$ of the universal enveloping algebra $U(\mathrm{sl}_2(K))$ of the 3 dimensional simple Lie algebra $\mathrm{sl}_2(K)$, where $C$ is the standard Casimir element and $0\neq \lambda \in K$. The structure of these algebras depend on $\lambda$ and $U(\mathrm{sl}_2(K))/(C-\lambda)$ is simple if $\lambda\neq n^2+n$ for any natural $n$. If $\lambda=n^2+n$ then $U(\mathrm{sl}_2(K))/(C-\lambda)$ has a unique non-trivial ideal of finite codimension $(n+1)^2$.

It is natural to ask if an analogue of the Dixmier result is true for quotients  $P(\mathrm{sl}_2(K))/(C-\lambda)$ of the Poisson enveloping algebra $P(\mathrm{sl}_2(K))$ of the simple Lie algebra $\mathrm{sl}_2(K)$, where $C$ is the standard Casimir element of $\mathrm{sl}_2(K)$ in $P(\mathrm{sl}_2(K))$.

Let  $L$  be an arbitrary  algebra with a skew-symmetric bilinear operation $[\,,]$ over a field $K$ of characteristic $0$ and let
$e_1,e_2\ldots$ be a linear basis of $L$. Then there exists a unique bracket $\{\,,\}$ satisfying the Leibniz identity on the
 polynomial algebra
$K[e_1, e_2,\ldots]$  such that $\{e_i,e_j\}=[e_i,e_j]$. The
algebra $\langle K[e_1,\ldots, ],\cdot, \{\,,\}\rangle$ will be called the {\em Poisson enveloping algebra} of $L$ and will be denoted by $P(L)$.

The algebra $P(L)$ satisfies the definition of generic Poisson algebras given in \cite{Shest}.
If $L$ is a Lie algebra then it is well known that $P(L)$ is a Poisson algebra. But if $M$ is a Malcev algebra then, in general, $P(M)$ is not a Malcev algebra with respect to  $\{\,,\}$.

The Poisson enveloping algebra $P(\mathrm{sl}_2(K))$ can also be obtained by the following construction. Let $C$ be an arbitrary non-constant element of the polynomial algebra $K[x,y,z]$.
 There exists a unique Poisson structure on $K[x,y,z]$ such that
\bes
\{x,y\}=\frac{\partial C}{\partial z}, \ \ \{y,z\}=\frac{\partial
C}{\partial x}, \ \ \{z,x\}=\frac{\partial C}{\partial y}.
\ees
We consider direct generalizations of these algebras for polynomial algebras in more than three variables. The new algebras are not Poisson algebras but $n$-Lie-Poisson algebras on the polynomial algebra $K[x_1,\ldots,x_{n+1}]$ ($n\geq 2$)
defined by
\bes
\{a_1,\ldots,a_n\}=J(a_1,\ldots,a_n,C)
\ees
 where $C$ is a fixed element of $K[x_1,\ldots,x_n]$ and $J$ is the Jacobian. Denote this $n$-Lie-Poisson algebra by $P_C$. Fortunately, this construction includes the $n$-Lie-Poisson enveloping algebras of all simple finite dimensional $n$-Lie algebras. Recall that all finite dimensional simple $n$-Lie algebras for $n>2$ have dimension $n+1$ \cite{Ling}.

 We describe the center of the $n$-Lie-Poisson algebra  $P_C$ and prove that the quotient algebra $P_C/(C-\lambda)$, where $(C-\lambda)$ is the ideal generated by $C-\lambda$ and $0\neq \lambda \in K$, is a central simple $n$-Lie-Poisson algebra if $C$ is a homogeneous element that is not a proper power of any nonzero element. We say that $C$ is not a proper power if $C$ cannot be represented in the form $C=\a c^k$ with $\a\in K$, $c\notin K$, and $k\geq 2$.

As we noticed above, the class of $n$-Lie-Poisson algebras of the form $P_C$ includes the Poisson enveloping algebra of $\mathrm{sl}_2(K)$ and the $n$-Lie-Poisson enveloping algebras of the simple $n$-Lie algebras of dimension $n+1$. We also prove similar results for the quotients of the Poisson enveloping algebra $P(\mathbb{M})$ of the simple exceptional Malcev algebra $\mathbb{M}$ of dimension $7$.

The paper is organized as follows. In Section 2 we give all necessary definitions and examples of Poisson and $n$-Lie Poisson algebras. In particular, the $n$-Lie-Poisson algebras of the form $P_C$ are defined. In Section 3 we describe the center of $P_C$. The class of algebras $S_C$ is defined in Section 4. In Section 5 we prove that $C-\lambda$ is irreducible if $C$ is closed and homogeneous and $0\neq \lambda\in K$. One important grading of $S_C$ is given in Section 6. In Section 7 we prove that $S_C$ is central simple. The standard central element $c_\mathbb{M}$ of the Poisson enveloping algebra $P(\mathbb{M})$ of the simple exceptional Malcev algebra $\mathbb{M}$ of dimension $7$ that generates the center of $P(\mathbb{M})$ is described in Section 8. Section 9 is devoted to the proof that $P(\mathbb{M})/(c_\mathbb{M}-\lambda)$ is central and simple.

\section{Definitions and examples}

\hspace*{\parindent}

A vector space $A$ over a field $K$ endowed with two bilinear
operations $x\cdot y$ (a multiplication) and $\{x,y\}$ (a Poisson
bracket) is called {\em a Poisson algebra} if $A$ is a commutative
associative algebra under $x\cdot y$, $A$ is a Lie algebra under
$\{x,y\}$, and $A$ satisfies the following identity (the Leibniz
identity):
\bes
\{x, y\cdot z\}=\{x,y\}\cdot z + y\cdot \{x,z\}.
\ees

There are two important classes of Poisson algebras.

{\em Example 1}. Symplectic Poisson algebras $P_n$. For each $n$ the algebra $P_n$ is a polynomial algebra
$K[x_1,y_1, \ldots,x_n,y_n]$ endowed with the Poisson bracket
defined by
\bes \{x_i,y_j\}=\delta_{ij}, \ \ \{x_i,x_j\}=0, \ \
\{y_i,y_j\}=0, \ees
where $\delta_{ij}$ is the Kronecker symbol
and $1\leq i,j\leq n$.

{\em Example 2}. Let  $L$  be an arbitrary  algebra with a skew-symmetric bilinear operation $[\,,]$ over a field $K$ of characteristic $0$ and let
$e_1,e_2\ldots$ be a linear basis of $L$. Then there exists a unique bracket $\{\,,\}$ satisfying the Leibniz identity on the
 polynomial algebra
$K[e_1, e_2,\ldots]$  such that $\{e_i,e_j\}=[e_i,e_j]$. The
algebra $\langle K[e_1,\ldots, ],\cdot, \{\,,\}\rangle$ will be called the {\em Poisson enveloping algebra} of $L$ and will be denoted by $P(L)$. Note that the bracket  $\{\,,\}$ of the algebra $P(L)$ depends on the structure of $L$ but
does not depend on a chosen basis.

This is a well known construction and if $L$ is a Lie algebra then $P(L)$ is a Poisson algebra.
But if $M$ is a Malcev algebra then, in general, $P(M)$ is not a Malcev algebra with respect to  $\{\,,\}$.

An associative and commutative algebra with a skew-symmetric bracket is called {\em generic Poisson} in \cite{Shest} if it satisfies only the Leibniz identity. Thus $P(L)$ is a generic Poisson algebra for any algebra $L$ with a skew-symmetric bilinear operation.

From now on we assume that $K$ is a field of characteristic zero and let $K[x,y,z]$ be a polynomial algebra in three variables.
Let $C$ be an arbitrary non-constant element of $K[x,y,z]$.

{\em Example 3}. There exists a unique Poisson structure on $K[x,y,z]$ such that
\bes
\{x,y\}=\frac{\partial C}{\partial z}, \ \ \{y,z\}=\frac{\partial
C}{\partial x}, \ \ \{z,x\}=\frac{\partial C}{\partial y}.
\ees

This example gives rise to another definition of the Poisson enveloping algebra $P(\mathrm{sl}_2(K))$ of the simple three dimensional Lie algebra $\mathrm{sl}_2(K)$.

{\em Example 4}. Let
\bes
h=
 \left[\begin{array}{cc}
 1  &  0  \\
 0  &  -1  \\
\end{array}\right], e=\left[\begin{array}{cc}
 0  &  1  \\
 0  &  0  \\
\end{array}\right], f=\left[\begin{array}{cc}
 0  &  0  \\
 1  &  0  \\
\end{array}\right],
\ees
be the standard basis of $\mathrm{sl}_2(K)$. We have
\bes
 \{h,e\}=2e, \  \{h,f\}=-2f, \ \{e,f\}=h.
\ees
Put $C=1/2 h^2+2ef$. It is easy to see that the Poisson bracket defined on $K[e,f,h]$ by $C$, as in Example 3, is the Poisson enveloping algebra of the Lie algebra $\mathrm{sl}_2(K)$.

Example 3 is a particular case of the following more general Poisson structure on polynomial algebras \cite{OR}.

{\em Example 5}. Let $h,C_1,\ldots,C_{n-2}$ be any fixed polynomials from  the polynomial algebra $K[x_1,\ldots, x_n]$. For every $f,g\in K[x_1,\ldots, x_n]$ set
\bes
\{f,g\}= h J(f,g,C_1,\ldots,C_{n-2}),
\ees
where $J$ is the Jacobian, i.e.,  $J(f_1,\ldots,f_n)=\mathrm{det}[\frac{\partial f_i}{\partial x_j}]$ for all $f_1\ldots,f_n\in K[x_1,\ldots, x_n]$.
The polynomial algebra $K[x_1,\ldots, x_n]$ with this bracket is a Poisson algebra.

This operation is a particular case of the following $n-m$-ary Nambu operations on $K[x_1,\ldots, x_n]$.

{\em Example 6}. Let $h,C_1,\ldots,C_m\in K[x_1,\ldots, x_n]$. For every $f_1,\ldots,f_{n-m}\in K[x_1,\ldots, x_n]$ put
\bes
\{f_1,\ldots,f_{n-m}\}= h J(f_1,\ldots,f_{n-m},C_1,\ldots,C_m).
\ees
It is known that this bracket defines on $K[x_1,\ldots, x_n]$ the structure of an $(n-m)$-Lie-Poisson algebra.

 Recall that an $n$-ary multilinear operation $[\cdot,\ldots,\cdot]$ on a vector space $V$  is called {\em skew-symmetric} if $[v_1,\ldots, v_n]=0$
whenever $v_i=v_j$ for $i\neq j$. Over fields of characteristic $\neq 2$, this is equivalent to
\bes
[v_{\s(1)}, \ldots,v_{\s(n)}] = \mathrm{sgn}(\s) [v_1,\ldots,v_n],
\ees
for any $\s$ from the symmetric group $S_n$, where $\mathrm{sgn}(\s)$
is the sign of the permutation $\s$.

A vector space $A$ over $K$ endowed with an $n$-ary skew-symmetric operation  $[\cdot,\ldots,\cdot]$ is called an {\em $n$-Lie} algebra if the following Jacobi identity holds:
\bee\label{f1}
[[u_1,\ldots, u_n], v_1,\ldots, v_{n-1}] = \sum_{i=1}^n
[u_1,\ldots, u_{i-1}, [u_i, v_1,\ldots, v_{n-1}], u_{i+1},\ldots, u_n],
\eee
for any $u_1,\ldots, u_n, v_1,\ldots, v_{n-1}\in V$.

These algebras were first introduced by V. Filippov \cite{Filippov}. Note that $2$-Lie algebras are Lie algebras.
It turns out that every finite dimensional simple $n$-Lie algebra over $K$ for $n>2$ is of dimension $n+1$ \cite{Ling}.

{\em Example 7}. Simple $n$-Lie algebras of dimension $n+1$. Let $V$ be a vector space with a linear basis $e_1,\ldots,e_{n+1}$. Put
\bes
[e_1,\ldots,\widehat{e_i},\ldots,e_{n+1}]=\a_i e_i, \ \ 1\leq i\leq n+1,
\ees
where $\widehat{e_i}$ means that $e_i$ is absent. The vector space $V$ with this bracket is an $n$-Lie algebra.
We denote this algebra by $\mathbb{L}_n(f)$, where
\bes
f=\frac{1}{2}\sum_{i=1}^{n+1} (-1)^{n-i+1}\a_i x_i^2\in K[x_1,\ldots, x_{n+1}].
\ees
More generally, if $f\in K[x_1,\ldots, x_{n+1}]$ is an arbitrary quadratic form, then we denote by $\mathbb{L}_n(f)$ the $n$-Lie algebra defined by
\bes
[e_1,\ldots,\widehat{e_i},\ldots,e_{n+1}]=(-1)^{n-i+1}\frac{\partial f}{\partial x_i}(e_1,\ldots,e_{n+1}),  \ \ 1\leq i\leq n+1.
\ees

It is known that $\mathbb{L}_n(f)$ is simple if and only if $f$ is non-degenerate \cite{Ling}. If $K$ is algebraically closed then
the simple $n$-Lie algebra $\mathbb{L}_n(f)$ of dimension $n+1$ is unique.

A vector space $P$ over $K$ endowed with a bilinear
operation $x\cdot y$ (a multiplication) and an $n$-ary multilinear operation $\{x_1,\ldots,x_n\}$ (an $n$-Lie-Poisson
bracket) is called an {\em $n$-Lie-Poisson} (see, for example \cite{Dzhumadildaev}) algebra if $P$ is a commutative
associative algebra under $x\cdot y$,
$P$ is a $n$-Lie algebra under
$\{x_1,\ldots,x_n\}$, and $P$ satisfies the following identity (the Leibniz
identity):
\bee\label{f2}
\{x_1,\ldots,x_{i-1}, x\cdot y,x_{i+1},\ldots,x_n\}=x\cdot\{x_1,\ldots,x_{i-1}, y,x_{i+1},\ldots,x_n\}\\
\nonumber  + \{x_1,\ldots,x_{i-1},x,x_{i+1},\ldots,x_n\}\cdot y.
\eee
The identity (\ref{f2}) means that the $n$-ary  operation $\{x_1,\ldots,x_n\}$ is a derivation with respect to the $i$th component.

{\em Example 8}. The operation
\bes
\{f_1,\ldots,f_n\}=J(f_1,\ldots,f_n)
\ees
defines the structure of an $n$-Lie-Poisson algebra on $K[x_1,\ldots, x_n]$. This is a simple $n$-Lie Poisson algebra and is an exact analogue of the Poisson symplectic algebra $P_1$ in the $n$-ary case when $n\geq 3$.

{\em Example 9}.  Poisson enveloping algebra $P(\frak{g})$ of an $n$-Lie algebra $\frak{g}$. Let $\frak{g}$ be an $n$-Lie algebra with a linear
basis $e_1,e_2,\ldots,e_k,\ldots$. Then there exists a unique $n$-ary skew-symmetric bracket satisfying (\ref{f2}) on the  polynomial algebra $K[e_1,
e_2,\ldots,e_k,\ldots]$ defined by
\bes
\{e_{i_1},\ldots,e_{i_n}\}=[e_{i_1},\ldots,e_{i_n}]
\ees
for all $i_1,\ldots,i_n$, where $[\cdot,\ldots,\cdot]$
is the bracket of the $n$-Lie algebra $\frak{g}$. The polynomial algebra $K[e_1,
e_2,\ldots,e_k,\ldots]$ with this bracket  will be denoted by $P(\frak{g})$ and is called the  {\em Poisson (or $n$-Lie Poisson) enveloping} algebra of $\frak{g}$.

Unfortunately, in general, $P(\frak{g})$ is not an $n$-Lie-Poisson algebra if $\frak{g}$ is an $n$-Lie algebra. But the next example  shows that $P(\mathbb{L}_n(f))$ is an $n$-Lie-Poisson algebra anyway.

We now define an exact extension of Example 3.

{\em Example 10}. Let $C\in K[x_1,\ldots, x_{n+1}]$ be an arbitrary polynomial. For every $f_1,\ldots,f_n\in K[x_1,\ldots, x_{n+1}]$ put
\bes
\{f_1,\ldots,f_n\}=  J(f_1,\ldots,f_n,C).
\ees
It follows that
\bes
\{x_1,\ldots,\widehat{x_i},\ldots,x_{n+1}\}=(-1)^{n-i+1}\frac{\partial C}{\partial x_i}, \ \ 1\leq i\leq n+1,
\ees
where $\widehat{x_i}$ means that $x_i$ is absent. The polynomial algebra $K[x_1,\ldots, x_{n+1}]$ with this $n$-ary bracket is
 an $n$-Lie-Poisson algebra (see, for example \cite{OR}). We denote this algebra by $P_C[x_1,\ldots, x_{n+1}]$.

 Notice that the Poisson enveloping algebra $P(\mathbb{L}_n(f))$ of the $n$-Lie algebra $\mathbb{L}_n(f)$ is isomorphic to $P_f[x_1,\ldots, x_{n+1}]$.
 Consequently, $P(\mathbb{L}_n(f))$ is an $n$-Lie-Poisson algebra.

An $n$-Lie-Poisson algebra $P$ is called a {\em strong $n$-Lie-Poisson} algebra  \cite{Dzhumadildaev} if it satisfies the following identity:
\bee\label{f3}
\sum_{i=1}^{n+1}
(-1)^i \{u_1,\ldots, u_{n-1}, v_i\} ¿ \{v_1,\ldots, \widehat{v_i}, \ldots,v_{n+1}\}=0.
\eee

\begin{lm}\label{l1}
$P_C[x_1,\ldots, x_{n+1}]$ is a strong $n$-Lie-Poisson algebra.
\end{lm}
\Proof
Let
\bes
S(u_1,\ldots,u_{n-1},v_1,\ldots,v_{n+1})=\sum_{i=1}^n
(-1)^i \{u_1,\ldots, u_{n-1}, v_i\} ¿ \{v_1,\ldots, \widehat{v_i}, \ldots,v_{n+1}\}.
\ees
Direct calculations show that $S(u_1,\ldots,u_{n-1},v_1,\ldots,v_{n+1})$  is a derivation with respect to each component $u_1,\ldots,u_{n-1},v_1,\ldots,v_{n+1}$. Consequently, in order to check the identity (\ref{f3}),
it is sufficient to check that
\bes
S(u_1,\ldots,u_{n-1},v_1,\ldots,v_{n+1})=0
\ees
 whenever $u_1,\ldots,u_{n-1},v_1,\ldots,v_{n+1}\in \{x_1,\ldots, x_{n+1}\}$. Since $S$ is skew-symmetric in the variables $v_1,\ldots,v_{n+1}$, we may assume that $v_i=x_i$ for all $1\leq i\leq n+1$. Since $S$ is skew-symmetric in the variables $u_1,\ldots,u_{n-2}$, without loss of generality, we may also assume that $u_i=x_i$ for all $1\leq i\leq n-1$. Then,
\bes
S(x_1,\ldots,x_{n-1},x_1,\ldots,x_{n+1})=(-1)^n \{x_1,\ldots, x_{n-1}, x_{n}\} ¿ \{x_1,\ldots, x_{n-1},x_{n+1}\}\\
 +(-1)^{n+1} \{x_1,\ldots, x_{n-1}, x_{n+1}\} ¿ \{x_1,\ldots, x_{n-1},x_{n}\}=0.
\ees
This proves the lemma. $\Box$

\section{The Center of $P_C$}

\hspace*{\parindent}

The {\em center}  $Z(P)$ of an $n$-Lie-Poisson algebra $P$ is the
set of all elements $c\in P$ such that $\{c,f_2,\ldots,f_n\}=0$ for all $f_2,\ldots,f_n\in
P$.  The elements of $Z(P)$ are called {\em Casimir} elements.

A polynomial $g\in K[x_1,\ldots, x_n]$ is called {\em closed} \cite{NN} if $K[g]$ is integrally closed in $K[x_1,\ldots, x_n]$. In fact, if $g$ is closed then $K[g]$ is algebraically closed in $K[x_1,\ldots, x_n]$. If $f$ and $g$ are algebraically dependent then \cite{Zaks} there exists $c$ such that $f,g\in K[c]$. If $g$ is closed then $K[g]\subseteq K[c]$ implies that $K[g]= K[c]$. Consequently, $f\in K[g]$.

A polynomial $g\in K[x_1,\ldots, x_n]$ is called a {\em root} of $f\in K[x_1,\ldots, x_n]$ if $f\in K[g]$. A root of minimal degree is called a {\em minimal root}.
For any $f\in K[x_1,\ldots, x_n]$ there exists \cite{NN} a closed polynomial $g\in K[x_1,\ldots, x_n]$ such that $f\in K[g]$. In this case $g$ is a minimal root of $f$. The minimal root $g$ is defined uniquely up to an affine transformation, i.e., up to an element of the form $\a g+\b, \a,\b\in K, \a\neq 0$. Consequently, the subalgebra $K[g]$ is uniquely defined.

\begin{pr}\label{pr1} Let $c$ be a minimal root of $C$. Then
$Z(P_C[x_1,\ldots, x_{n+1}])= K[c]$.
\end{pr}
\Proof We follow the proof of Lemma 5 from \cite{MLTU-EP}. Put $P=P_C[x_1,\ldots, x_{n+1}]$. It is easy to check that $C\in Z(P)$.
Consequently, $c\in Z(P)$.  Assume that $Z(P)$ contains an element $f$ such that $c$ and $f$ are algebraically independent.
Since the transcendence degree of $K[c, f]$ over $K$ is $2$ and the transcendence degree of $P$ over $K$ is $n+1$,
we can find $n-1$ elements $g_1,\ldots,g_{n-1}$ of $P$  such that $c,f,g_1,\ldots,g_{n-1}$ are algebraically independent over $K$.
We get
\bes
 \{g_1,\ldots,g_{n-1}, f\}=0
\ees
since $f$ is in the center of $P$. By the definition of the bracket in $P$ we also get
\bes
 J(g_1,\ldots,g_{n-1}, f, C)=0
\ees
in $K[x_1,\ldots, x_{n+1}]$. Consequently,
\bes
 J(g_1,\ldots,g_{n-1}, f, c)=0.
\ees
This implies (see, for example \cite{SU04}) that $g_1,\ldots,g_{n-1}, f, c$ are algebraically dependent over $K$.

This contradiction proves that if $f\in Z(P)$ then $f,c$ are algebraically dependent over $K$. Since $c$ is closed it follows that $f\in K[c]$. $\Box$

\section{Algebras $S_C$}

\hspace*{\parindent}

In Sections 4--7 we assume that $C\in K[x_1,\ldots, x_{n+1}]$ is a closed homogeneous polynomial of degree $m\geq 1$. The well known Euler formula states that
\bes
\sum_{i=1}^{n+1} x_i\frac{\partial C}{\partial x_i}=mC.
\ees
If $C$ is homogeneous then $C$ is closed if and only if $C$ cannot be represented in the form $C=\a c^t$, where $\a\in K$, $c\in K[x_1,\ldots, x_{n+1}]$, and $t\geq 2$.

Denote by $S_C[x_1,\ldots, x_{n+1}]$ the quotient algebra of $P_C[x_1,\ldots, x_{n+1}]$ by the Poisson ideal $(C-\lambda)$ generated by $C-\lambda$ where $0\neq \lambda\in K$. Since $C$ belongs to the center of $P_C[x_1,\ldots, x_{n+1}]$, it follows that
$S_C[x_1,\ldots, x_{n+1}]$ is a Poisson algebra on the commutative and associative algebra $K[x_1,\ldots, x_{n+1}]/(C-\lambda)$, where $(C-\lambda)$ is the ideal of $K[x_1,\ldots, x_{n+1}]$ generated by $C-\lambda$.

The main result of the following three sections is to prove that  $S_C[x_1,\ldots, x_{n+1}]$ is a central simple $n$-Lie-Poisson algebra if $C$ is a closed homogeneous polynomial and $\lambda\neq 0$.

{\em Example 11}. If $C\in K[x,y,z]$ is linear then we can make $C=z$. Consequently,
\bes
\{x,y\}=1, \ \ \{y,z\}=0, \ \ \{z,x\}=0,
\ees
in $P_C[x,y,z]$. Then $S_C[x,y,z]$ is isomorphic to the symplectic Poisson algebra $P_1$ and it is a simple Poisson algebra.

Notice that if $C=c^t,t\geq 2$ and $\lambda=1$, then $S_C$ cannot be simple since $(c-1)$ induces a non-zero proper ideal of $S_C$.
In fact, if $C-\lambda$ is reducible then $S_C$ cannot be simple.

\section{Irreducibility of $C-\lambda$}

\hspace*{\parindent}

\begin{pr}\label{pr2}
Let $C\in K[x_1,\ldots, x_{n+1}]$ be a closed homogeneous non-constant polynomial. Then $C-\lambda$,  where $0\neq \lambda\in K$, is irreducible.
\end{pr}
\Proof Notice that if $C$ is closed in $K[x_1,\ldots, x_{n+1}]$ then $C$ is still closed in $\overline{K}[x_1,\ldots, x_{n+1}]$, where $\overline{K}$ is the algebraic closure of $K$. Indeed, suppose that $C=c^k$ for some $k\geq 2$ and $c\in \overline{K}[x_1,\ldots, x_{n+1}]$. Applying some linear transformations, if necessary,
we can assume that $C=\a C'$, where $C'$ can be written as an element of $K[x_1,\ldots,x_n][x_{n+1}]$ with the leading coefficient 1. Put also $c=\b c'$, where $\b$ is the leading coefficient of $c$. Then $C=c^k$ implies that $\a=\b^k$ and $C'=c'^k$. Assume that
\bes
C'= \sum_{i=0}^{mk}a_i x_{n+1}^i, c'= \sum_{i=0}^{m} b_i x_{n+1}^i , a_{mk}=b_m=1.
\ees
Since $C'=c'^k$ it follows that
\bes
a_s=\sum_{i_1+\ldots +i_k=s} b_{i_1}\ldots b_{i_k}.
\ees
Suppose that $s<mk$ and $s=m(k-1)+r$, where $0\leq r<m$. Then we can write
\bee\label{f4}
a_s=S+kb_r
\eee
where $S$ is the sum of products of the form $b_{i_1}\ldots b_{i_k}$  with $i_1+\ldots +i_k=s$ and at least two of $i_1,\ldots,i_k$ are different from $m$.
We show that $i_1,\ldots,i_k>r$. In fact, if at least one of them is less than or equal to $r$ and one more is less than or equal to $m-1$, then we get
$i_1+\ldots +i_k\leq r+(m-1)+m(k-2)= m(k-1)+r-1<s$.

Consequently, $S$ does not contain coefficients $b_i$ with $i<r$ in (\ref{f4}). Recall that $a_i\in K[x_1,\ldots,x_n]$ for all $i$ and $b_m=1$. Then (\ref{f4}) allows us to prove that $b_r\in K[x_1,\ldots,x_n]$ by inverse induction on $r$. Therefore $c'\in K[x_1,\ldots, x_{n+1}]$ and
$C=\a c'^k$ is not closed.

Hence, we may assume that $K$ is an algebraically closed field of characteristic zero. We can also assume that $\lambda=1$ by changing $C$ to $C/\lambda$.
 Suppose that
\bes
C-1=fg,  \ \ \ \ \deg\,f, \deg\,g\geq 1.
\ees
The substitution $x_i\mapsto x_i/T$, where $T$ is a new variable, leads to a nontrivial homogeneous decomposition
\bee\label{f5}
T^m-C=FG,  \ \ \ \ \deg\,F, \deg\,G\geq 1,
\eee
where $m$ is the degree of $C$. Consider (\ref{f5}) as a decomposition of the polynomial $T^m-C$ in one variable $T$ over the quotient field  $K(x_1,\ldots,x_{n+1})$. We may assume that the leading coefficients of $F$ and $G$ are $1$'s.

Let $Q$ be any extension of the field $K(x_1,\ldots,x_{n+1})$ containing $b$ such that $b^m=C$. Then
\bes
T^m-C=(T-\varepsilon_1 b)\ldots (T-\varepsilon_m b),
\ees
where $\varepsilon_1,\ldots,\varepsilon_m$ are all root of unity of degree $m$. Without loss of generality we may assume that
\bee\label{f6}
F=(T-\varepsilon_1 b)\ldots (T-\varepsilon_k b), \ \ G=(T-\varepsilon_{k+1} b)\ldots (T-\varepsilon_m b).
\eee
Recall that (\ref{f5}) is a homogeneous decomposition in $K[x_1,\ldots, x_{n+1}][T]$ and all coefficients of $F$ and $G$, as a polynomial in $T$, belong to $K[x_1,\ldots, x_{n+1}]$. Then (\ref{f6}) implies that all coefficients of $F$ and $G$ have the form $\lambda b^k$ for some $\lambda\in K$ and a nonnegative integer $k$. Consequently, any two coefficients of $F$ and $G$ are algebraically dependent. By \cite{Zaks}, there exists an element $c\in K[x_1,\ldots, x_{n+1}]$ such that $K[c]$ contains all coefficients of $F$ and $G$. Obviously, we can choose homogeneous $c$ since all coefficients of  $F$ and $G$ are homogeneous. Then we  get $C\in K[c]$ and $\deg\,C>\deg\,F\geq \deg\,c$. This means that $C$ is not closed. $\Box$

\section{A grading of $S_C$}

\hspace*{\parindent}

As before, we assume that $C\in K[x_1,\ldots, x_{n+1}]$ is a closed homogeneous polynomial of degree $m\geq 1$. We have
\bes
K[x_1,\ldots, x_{n+1}]=\bigoplus_{r=0}^{m-1} A_r,
\ees
where $A_r$ is the linear span of all homogeneous polynomials of degree $k$ such that $k\equiv r (\mathrm{mod}\, m)$ for all $0\leq r\leq m-1$.
Obviously, $A_rA_{r'}\subseteq A_{r\oplus r'}$, where $r\oplus r'$ denotes addition in $\mathbb{Z}_m=\{0,1,\ldots,m-1\}$, i.e., this decomposition makes
 $K[x_1,\ldots, x_{n+1}]$ a $\mathbb{Z}_m$-graded algebra.
The elements of $A_r$ will be called {\em $m$-homogeneous}.

Since  $C-\lambda$ is $m$-homogeneous, this grading induces the $\mathbb{Z}_m$-grading of the quotient algebra
\bes
K[x_1,\ldots, x_n]/(C-\lambda)=\bigoplus_{r=0}^{m-1} B_r,
\ees
where $B_r$ is the image of $A_r$. The elements of $B_r$ will also be called {\em $m$-homogeneous}.

For any $f\in K[x_1,\ldots, x_{n+1}]$ denote by $\overline{f}$ the image of $f$ in $K[x_1,\ldots, x_{n+1}]/(C-\lambda)$.

\begin{lm}\label{l2} Let $f\in A_r$ and $\overline{f}\neq 0$. Then there exists a homogeneous and $m$-homogeneous polynomial $f'\in A_r$ such that $\overline{f}=\overline{f'}$ and $\deg\,f=\deg\,f'$.
\end{lm}
\Proof We have
\bes
f=f_{n_k}+f_{n_{k-1}}+\ldots +f_{n_0},
\ees
where $f_{n_i}$ is homogeneous of degree $n_i$, $n_i\equiv r (\mathrm{mod} m)$ for all $i$, and $n_k>n_{k-1}>\ldots > n_0$. Set
\bes
f'=f_{n_k}+f_{n_{k-1}}(C/\lambda)^{t_{k-1}}+\ldots +f_{n_0}(C/\lambda)^{t_0},
\ees
where $n_k-n_i=mt_i$ for all $0\leq i\leq k-1$. Then $f'$ satisfies the required conditions. $\Box$

\begin{lm}\label{l3} The decomposition
\bee\label{f7}
S_C=\bigoplus_{r=0}^{m-1} B_r
\eee
defines a grading of the $n$-Lie-Poisson algebra $S_C$.
\end{lm}
\Proof If $f_i\in A_{r_i}$ for all $1\leq i\leq n$ then
\bes
\{\overline{f_1},\ldots,\overline{f_n}\}=\overline{\{f_1,\ldots,f_n,C\}}\in B_r,
\ees
where $(r_1-1)+\ldots+(r_n-1)+(m -1)\equiv r(\mathrm{mod}\, m)$ and $0\leq r\leq m-1$. $\Box$

\section{Simplicity of $S_C$}

\hspace*{\parindent}

Let $P$ be an $n$-Lie-Poisson algebra over $K$. We say that $P$ is {\em central} over $K$ if $Z(P)=K$.

\begin{theor}\label{t1} Let $C\in K[x_1,\ldots, x_{n+1}]$ be an arbitrary non-constant homogeneous closed polynomial.
Then the $n$-Lie-Poisson algebra $S_C[x_1,\ldots, x_{n+1}]$ is a central simple $n$-Lie-Poisson algebra over $K$.
\end{theor}
\Proof Since $C$ belongs to the center of $P_C=P_C[x_1,\ldots, x_{n+1}]$, it follows that
$S_C=S_C[x_1,\ldots, x_{n+1}]$ is a Poisson algebra on the commutative and associative algebra\\
 $K[x_1,\ldots, x_{n+1}]/(C-\lambda)$, where $(C-\lambda)$ is the ideal of $K[x_1,\ldots, x_{n+1}]$ generated by $C-\lambda$.
By Proposition \ref{pr2}, $C-\lambda$ is irreducible and $(C-\lambda)$ is a prime ideal.
The Krull dimension theory \cite[Chapter 5]{Matsumura} says that the height $\mathrm{hgt}(C-\lambda)=1$ since $(C-\lambda)\supsetneqq (0)$.

Assume that $S_C$ is not simple. Then $(S_C, \cdot)$ contains a
nonzero ideal $I$ such that $\{S_C,I\}\subseteq I$. Consequently,
for every derivation $\mathrm{ad}_{u_1,\ldots,u_{n-1}} : a\mapsto \{u_1,\ldots,u_{n-1},a\}$ of $S_C$ we get
$\mathrm{ad}_{u_1,\ldots,u_{n-1}}(I)\subseteq I$. The set of the ideals with these properties contains a maximal
ideal $J$. It is known \cite{Posner60} that $J$ is a prime ideal of $(S_C, \cdot)$.

Let $\widehat{J}$ be the pre-image of $J$ in $K[x_1,\ldots, x_{n+1}]$. Then
$\widehat{J}$ is a prime ideal of $K[x_1,\ldots, x_{n+1}]$ and we have
$\widehat{J} \supsetneq (C-1)\supsetneqq (0)$. Consequently,
$\mathrm{hgt}(\widehat{J})\geq 2$. This implies (see \cite[Chapter 5]{Matsumura}) that
\bes
\mathrm{tr.deg}(K[x_1,\ldots, x_{n+1}]/\widehat{J})\leq n-1.
\ees
Then $\mathrm{tr.deg}(S_C/J)\leq n-1$. Consequently, any $n$
elements of $S_C/J$ are algebraically dependent.

Let $g_1,\ldots,g_n$ be arbitrary elements of $S_C/J$. Then there exists a nonzero polynomial $F(y_1,\ldots,y_n)\in K[y_1,\ldots,y_n]$ such that
\bes
F(g_1,\ldots,g_n)=0.
\ees
Of course, we may assume that $F$ has the minimal possible degree in $y_n$. Then
\bes
0=\{g_1,\ldots,g_{n-1},F(g_1,\ldots,g_n)\}=\{g_1,\ldots,g_{n-1},g_n\}\frac{\partial F}{\partial y_n}(g_1,\ldots,g_n)
\ees
implies that
\bes
\{g_1,\ldots,g_{n-1},g_n\}=0
\ees
since $\frac{\partial F}{\partial y_{n-1}}(g_1,\ldots,g_n)\neq 0$.
Consequently, the
Poisson bracket on $S_C/J$ is trivial, i.e.,
\bes
\{x_1,\ldots,\widehat{x_i},\ldots,x_{n+1}\}=0
\ees
in $S_C/J$ for all $0 \leq i\leq n+1$. Therefore
\bes
\{x_1,\ldots,\widehat{x_i},\ldots,x_{n+1}\}=(-1)^{n-i+1}\frac{\partial C}{\partial x_i}\in J, \ \ 1\leq i\leq n+1.
\ees
By Euler's formula, we get
\bes
\sum_{i=1}^{n+1} x_i\frac{\partial C}{\partial x_i}=mC \in J,
\ees
where $m=\deg\,C$. Since $C-\lambda\in J$ it follows that $J=S_C$.

This proves that $S_C$ is simple.

Suppose that $g\in Z(S_C)$. Consider the grading (\ref{f7}) of $S_C$. Obviously, all $m$-homogeneous components of $g$ belong to the center of $S_C$ too.
Consequently, we may assume that $g$ is $m$-homogeneous and $g\in B_r$. By Lemma \ref{l2}, there exists a homogeneous and $m$-homogeneous element $f\in A_r$ such that $g=\overline{f}$. We get
\bes
\{x_1,\ldots,\widehat{x_i},\ldots,\widehat{x_j},\ldots,x_{n+1},g\}=\overline{\{x_1,\ldots,\widehat{x_i},\ldots,\widehat{x_j},\ldots,x_{n+1},f,C\}}=0
\ees
for all $i<j$. Notice that
\bes
\{x_1,\ldots,\widehat{x_i},\ldots,\widehat{x_j},\ldots,x_{n+1},f,C\}=\pm (\frac{\partial f}{\partial x_i}\frac{\partial C}{\partial x_j}-\frac{\partial f}{\partial x_j}\frac{\partial C}{\partial x_i}).
\ees
Hence,
\bes
\frac{\partial f}{\partial x_i}\frac{\partial C}{\partial x_j}-\frac{\partial f}{\partial x_j}\frac{\partial C}{\partial x_i}\in (C-\lambda).
\ees
But the ideal $(C-\lambda)$ cannot contain homogeneous elements except $0$. Therefore
\bes
\frac{\partial f}{\partial x_i}\frac{\partial C}{\partial x_j}-\frac{\partial f}{\partial x_j}\frac{\partial C}{\partial x_i}=0
\ees
in $K[x_1,\ldots, x_{n+1}]$ for all $1\leq i<j\leq n+1$. It is well known (see, for example \cite{SU04}) that this implies algebraic dependency of $f$ and $C$. Since $C$ is closed it follows that $f\in K[C]$. Then $g\in K$ in $S_C$.
$\Box$

Recall that $C=1/2 h^2+2ef$ is the Casimir element of the Poisson algebra $P(\mathrm{sl}_2(K))$ in the standard basis $e,f,h$ (see, Example 4).
\begin{co}\label{c1} Let $P(\mathrm{sl}_2(K))$ be the Poisson enveloping algebra of $\mathrm{sl}_2(K)$ with the standard Casimir element $C=1/2 h^2+2ef$. Then the quotient Poisson algebra $P(\mathrm{sl}_2(K))/(C-\lambda)$ is central and simple for any $0\neq \lambda\in K$.
\end{co}

This is a Poisson analogue of the Dixmier \cite{Dixmier73} theorem on the quotients $U(\mathrm{sl}_2(K))/(C-\lambda)$ of the universal enveloping algebra $U(\mathrm{sl}_2(K))$ mentioned in the Introduction.

{\em Example 12}. (Elliptic Poisson algebras) Let
$C=\frac{1}{3}(x^3+y^3+z^3)-\alpha xyz$, $\alpha \in K$. Then $P_C[x,y,z]$ is called an {\em elliptic Poisson} algebra and
the bracket is defined by
\bes \{x,y\}=-\a xy+z^2, \ \ \{y,z\}=-\a
yz +x^2, \ \ \{z,x\}=-\a zx+y^2.
\ees

Automorphisms of the elliptic Poisson algebras were described in \cite{MLTU-EP}.

 \begin{co}\label{c2} Let $P_C[x,y,z]$ be the elliptic Poisson algebra with the standard Casimir element $C=\frac{1}{3}(x^3+y^3+z^3)-\alpha xyz$.
 Then the quotient Poisson algebra $P_C[x,y,z]/(C-\lambda)$ is central and simple for any $0\neq \lambda\in K$.
\end{co}

Recall that every $n+1$-dimensional $n$-Lie algebra $\mathbb{L}(f)$ is defined by a quadratic form $f\in K[x_1,\ldots,x_{n+1}]$ (see, Example 7).
Moreover, $\mathbb{L}(f)$ is simple if and only if $f$ is non-degenerate. Obviously, non-degenerate quadratic polynomials are closed if $n\geq 2$.
\begin{co}\label{c2} Let $P(\mathbb{L}(f))$ be the $n$-Lie-Poisson enveloping algebra of the simple $n$-Lie algebra $\mathbb{L}(f)$. Then the quotient Poisson algebra $P(\mathbb{L}(f))/(f-\lambda)$ is central and simple for any $0\neq \lambda\in K$.
\end{co}

\section{The standard Casimir element and the center of $P(\mathbb{M})$}

\hspace*{\parindent}

In this section we describe the standard Casimir element of the Poisson enveloping algebra $P(\mathbb{M})$ of an  exceptional simple Malcev algebra $\mathbb{M}$ of dimension $7$. This element generates the center of $P(\mathbb{M})$.

Recall that an algebra $M$ over a field of characteristic $\neq 2,3$ endowed with a binary skew-symmetric operation $[\,,]$ is called a {\em Malcev}
algebra if it satisfies the identity
$$[J(x,y,z),x]=J(x,y,[x,z]),$$
where $J(x,y,z)=[[x,y],z]-[[x,z],y]-[x,[y,z]]$ is called the {\em Jacobian} of elements $x,y,z\in M$. (Do not confuse this Jacobian with the Jacobian of a system of polynomials used in the preceding sections.)

In 1955 A.I. Mal'cev \cite{Malcev} introduced Malcev algebras as the tangent algebras of locally analytic Moufang loops. Every Lie algebra is a Malcev algebra.
An important example of a Malcev algebra is the traceless subspace of the Cayley-Dickson algebra with respect to the commutator operation. In general, for any alternative algebra $A$ the commutator algebra $A^{(-)}$ is a Malcev algebra.

The (nonassociative) universal enveloping algebra $U(M)$ of an arbitrary Malcev algebra $M$ is defined in \cite{IzSh}. The algebra $U(M)$ is not alternative in general,  but it has a nonassociative Poincare-Birkhoff-Witt basis.

Let $\mathbb{O}$ be a Cayley-Dickson algebra over a field $K$ of characteristic $\neq 2$ and let $\mathrm{tr}$ be the trace function (see, for example \cite{KS,KBKA}) on
$\mathbb{O}$. Then $\mathbb{O}=K \oplus \mathbb{O}_0$, where $\mathbb{O}_0=\{x\in \mathbb{O}\,|\,\mathrm{tr}(x)=0\}$. The product of elements $a,b\in \mathbb{O}_0$ is defined by
 $$
 a\cdot b=(a,b)+[a,b],
 $$
where $(\,,)$ is a symmetric non-degenerate bilinear form on $\mathbb{O}_0$ and $[\,,]$ is the commutator on $\mathbb{O}_0$. The algebra $\mathbb{M}=(\mathbb O_0,[\,,])$ is a seven dimensional Malcev algebra and it is simple if the characteristic of $K$ is $\neq 2,3$. Moreover, every central simple Malcev algebra that is not a Lie algebra over a field of characteristic $\neq 2,3$ is isomorphic to an algebra of type $\mathbb M$.

If $K$ is algebraically closed then $\mathbb M$ has a canonical basis $e_1,\ldots, e_7$ with the multiplication table
 \begin{eqnarray}\label{f8}
[e_i,e_{i+1}]=e_{i+3}, \  \hbox{where} \ e_{7+j}=e_{j} \ \hbox{for any} \  j>0,
 \end{eqnarray}
that is invariant with respect to the cyclic permutation of triples of indexes. Then
$$
c_{\mathbb M}=e_1^2+\ldots+e_7^2
$$
belongs to the center (see the definition below) of the universal enveloping algebra $U(\mathbb M)$.

If $K$ is an algebraically closed field then one can also choose a so called {\em splittable} basis $h,x,y,z,x',y',z'$ of $\mathbb M$ with multiplication table

\bes
 &&[h,x]=2x,\ [h,y]=2y,\ [h,z]=2z, \\ \
 &&[h,x']=-2x',\ [h,y']=-2y',\ [h,z']=-2z',\\ \
 &&[x,x']=[y,y']=[z,z']=h,\\ \
 && [x,y]=2z',\ [y,z]=2x',\ [z,x]=2y',\\ \
 &&[x',y']=-2z,\ [y',z']=-2x,\ [z',x']=-2y,
\ees
where all absent products are zeroes \cite{K89}. Each of the triples
 $\{h,x,x'\},\ \{h,y,y'\},\ \{h,z,z'\}$ forms a splittable basis for $\mathrm{sl}_2(K)$.
 This splittable basis of $\mathbb M$ can be obtained from the canonical basis $e_i$ as follows:
\bes
 && x=ie_2+e_4,\ y=ie_3+e_7, \ z=e_5-ie_6, \\
 &&x'=ie_2-e_4,\ y'=ie_3-e_7,\ z'=-e_5-ie_6.
\ees
The central element $c_{\mathbb M}\in Z(U(\mathbb M))$ can be written as
\bee\label{f9}
c_{\mathbb M}=-(x\circ x'+y\circ y'+z\circ z'+\tfrac{1}{4}h^2),
\eee
where $a\circ b =\tfrac{1}{2}(ab+ba)$.

In the case of an arbitrary field $K$ the algebra $\mathbb M=\mathbb M(\alpha, \beta,\gamma)$
is defined by a triple of nonzero parameters
 $\alpha, \beta,\gamma\in K$ and has a basis $f_1,\ldots,
f_7$, which can be formally expressed as
\begin{eqnarray*}
&& f_{1}=\sqrt\alpha \,e_{1},\,f_{2}=\sqrt\beta
\,e_{2},\,f_{3}=\sqrt\gamma \,e_{3},\,
f_{4}=\sqrt{\alpha\beta}\, e_{4},\\
&&f_{5}=\sqrt{\beta\gamma}\, e_{5},\, f_{6}=\sqrt{\alpha\beta\gamma}
\,e_{6},\,f_{7}=\sqrt{\alpha\gamma}\, e_{7};
\end{eqnarray*}
where the elements $e_{i}$ are multiplied according to \ref{f8}. The table of multiplication for the elements
contains only positive integer powers of the parameters $\alpha,\beta,\gamma$:
$$
[f_{1},f_{2}]=f_{4},\ [f_{2},f_{4}]=\beta\,f_{1},\
[f_{5},f_{6}]=\beta\gamma\,f_{1},\ldots
$$

The central element of $U(\mathbb M)$ corresponding to $c_{\mathbb M}$ is
$$
c_{\mathbb M}(\alpha,\beta,\gamma)=\alpha \beta\gamma c_{\mathbb
M}=\beta\gamma f_1^2+\alpha\gamma f_2^2+ =\alpha \beta f_3^2+=\gamma
f_4^2+\alpha f_5^2+ f_6^2 +\beta f_7^2.
$$

Recall that the associative center $N(A)$, the commutative center $K(A)$, and
the center $Z(A)$ of an arbitrary algebra $A$ are defined as follows:
\begin{eqnarray*}
N(A)&=&\{n\in A|\,(n,a,b)=(a,n,b)=(a,b,n)=0 \ \text{for all} \  a,b\in
A\},\\
 K(A)&=&\{k\in A|\,[k,a]=0 \ \text{for all} \ a\in A\},\\
 Z(A)&=&N(A)\cap K(A).
\end{eqnarray*}

The center of the universal algebra $U(\mathbb M)$ of a seven dimensional Malcev algebra $\mathbb M$ is described in \cite{ZhelShest}. If $\mathbb M={\mathbb M}(\alpha,\beta,\gamma)$ then the associative center of the algebra $U(\mathbb M)$ coincides with its center and is a polynomial algebra in one variable $c_{\mathbb M}(\alpha,\beta,\gamma)$.

\begin{pr}\label{pr3}  Let $\mathbb M={\mathbb M}(\alpha,\beta,\gamma)$ be a seven dimensional simple Malcev algebra over a field $K$ of characteristic zero with nonzero parameters
$\alpha, \beta,\gamma\in K$. Then the center of the Poisson enveloping algebra $P(\mathbb M)$ is a polynomial algebra in one variable $c_{\mathbb M}(\alpha,\beta,\gamma)$.
\end{pr}
The proof of Theorem 2 from \cite{ZhelShest} contains the proof of Proposition \ref{pr3}.

\section{Simple quotients of $P(\mathbb{M})$}

\hspace*{\parindent}

Let $\langle A,\cdot, \{\,,\}\rangle$ be a generic Poisson algebra. As in the case of Poisson algebras, $A$ is
called {\em simple} if $A$ does not contain proper ideals $I$ such that $\{A,I\}\subseteq I$.

\begin{rem}\label{r1}
Let  $A$ be a unital generic Poisson algebra and let $F$ be an extension of the field $K$. If the generic Poisson algebra $A\otimes_KF$ is simple
then $A$ is simple too. If $A\otimes_KF$ is central over $F$ then $A$ is central over $K$.
\end{rem}
\Proof
Assume that $A\otimes_KF$ is simple and
let $I$ be a nonzero ideal of the algebra $A$. Then
$I\otimes_KF$ is a nonzero ideal of  $A\otimes_KF$. Consequently,
 $I\otimes_KF=A\otimes_KF$ and
\bes
1\otimes
1=\sum_{i=1}^nr_i\otimes f_i,
\ees
 where $r_i\in I, f_i\in F$. We
can assume that $f_1,\ldots, f_n$ are linearly independent over $K$.
Then  $1=\sum_i\alpha_if_i$, for some $\alpha_i\in K$. From here we get
\bes
\sum_{i=1}^n(r_i-\alpha_i\cdot 1)\otimes f_i=0.
\ees
 Consequently,
$r_i=\alpha_i\cdot 1$ for some index $i$. Therefore $I=A$.

If $a\in Z(A)\setminus K$ then $a\otimes 1\in Z(A\otimes_KF)\setminus F$. $\Box$

Let $ \lambda\in K$ and let $(c_\mathbb M-\lambda)$ be the ideal of $\langle P(\mathbb M), \cdot\rangle$
generated by $c_{\mathbb
M}-\lambda$. Since $c_{\mathbb M}-\lambda$ belongs to the center of $\langle P(\mathbb M),\{\,,\}\rangle$ it follows
that $(c_\mathbb M-\lambda)$ is an ideal of the generic Poisson algebra $\langle P(\mathbb M),\cdot,\{\,,\}\rangle$.  Denote by $S=P(\mathbb M)/(c_\mathbb M-\lambda)$ the quotient generic Poisson algebra.

\begin{lm} \label{l4} Let $K$ be a field of characteristic $0$ and $0\neq \lambda\in K$. Then the generic Poisson algebra
 $S=P(\mathbb M)/(c_\mathbb M-\lambda)$ is simple.
\end{lm}
\Proof By Remark \ref{r1} we may assume that $K$ is algebraically closed. Therefore $c_{\mathbb M}\in P(\mathbb M)$ can be written in the form (\ref{f9}).
The image of the splittable basis $x,x',y,y',z,z',h$ in $S$ will be denoted by the
same symbols.

Assume that $S$ is not simple.
 Then the algebra $\langle S,\cdot\rangle$ contains a proper ideal $I$ such that $\{S,I\}\subseteq I$.
For any $u\in S$ the mapping $d_u:a\mapsto \{u,a\}$ is a
derivation of the algebra $\langle S,\cdot\rangle$. Moreover, $d_u(I)\subseteq
I$ for any $u\in S$. Between such ideals we can choose a maximal
ideal $J$  since $\langle S,\cdot\rangle$ is Noetherian. It is known \cite{Posner60} that $J$ is a prime ideal of $\langle S,\cdot\rangle$.
 Notice that $\mathbb{M}\cap J=0$ since otherwise the simplicity of $\mathbb M$ implies $J=S$.

By Proposition \ref{pr1}, $c_{\mathbb M}-\lambda$ is irreducible
over $K$. Consequently, $S$ is an integral domain. Set
$A=K[x,x',y,y',z,z']$. Then  $A\cap J\neq 0$. Indeed, let
$0\neq r\in J$. Since $S=A+hA$, we can write $r=a+hb$, where
$a,b\in A$. Notice that $A\cap hA=0$. Consequently, $a-hb\neq 0$ and
\bes
0\neq r(a-hb)=a^2-h^2b^2\in A\cap J.
\ees
Hence, $A\cap J\neq 0$.

Let $J(a,b,c)$ be the Jacobian with respect to the bracket $\{\,, \}$. Then the mapping
$a \mapsto
J(a,b,c)$ is a derivation of $\langle S,\cdot\rangle$ for any fixed $b,c\in S$.

Suppose that $0\neq f\in K[x,x',y',z']\cap J$. Let
$f=\sum_{i=0}^nf_ix^i,$ where $f_i\in K[x',y',z']$.
We have
\bes
J(x,y,h)=12z',J(x',y,h)=J(y',y,h)=J(z',y,h)=0.
\ees
Using this, we get
\bes
J(f,y,h)=\sum_{i=0}^n12f_iix^{i-1}z'\in J.
\ees
Since $z'\not \in J$ it follows that $\sum_{i=0}^nf_iix^{i-1}\in J$. This allows us to reduce the degree of $f$ in $x$. Consequently, we may assume that
$f\in K[x',y',z']\cap
J$.

Let   $f=\sum_{i=0}^nf_ix'^i$, where  $f_i\in K[y',z']$.
Since $\{x,y'\}=\{x,z'\}=0 $ and  $\{x,x'\}=h$, direct calculation gives that
\bes
\{x,f\}=\sum_{i=0}^nf_iix'^{i-1}h\in J.
\ees
Since $h\not \in J$ it follows that $\sum_{i=0}^nf_iix'^{i-1}h\in J$. This allows us to reduce the degree of $f$ in $x'$.
Consequently, we may assume that $f\in K[y',z']\cap J$.

 Continuing similarly, we can get $f\in K$. Therefore, we may assume that $K[x,x',y',z']\cap J=0$.

Suppose that $f=\sum_{i=0}^nf_{i}y^i$, where $f_{i}\in K[x,x',y',z']$.
We have
\bes
J(x', x,h)=J(y', x, h)=J(z',x,h)=0, \ J(y,x,h)=-12z', \
J(z,x,h)=12y'.
\ees
 Then
 \bes
 J(f,x,h)=\sum_{i=0}^n-12f_{i}iy^{i-1}z'\in J.
 \ees
Since $z'\not \in J$ it follows that $\sum_{i=0}^nf_{i}iy^{i-1}\in J$. This allows us to reduce the degree of $f$ in $y$ if $f$ does not depend on $z$.
Similarly, we can reduce the degree of $f$ in $z$ if $f$ does not depend on $y$.
Consequently, we may assume that $f$ depends on both $y$ and $z$.

Let $0\neq f\in A\cap J$ and
\bes
f=\sum_{k=0}^n\sum_{i+j=k}f_{ij}y^iz^j,
\ees
where $f_{ij}\in K[x,x',y',z']$.
We choose $f$ with the minimal possible degree $\deg_yf+\deg_zf >0$. Since
$J(f_{ij},x,h)=0$ it follows that
\bes
J(f,x,h)=\sum_{k=1}^{n}\sum_{i+j=k}12f_{ij}(-iz'y^{i
-1}z^{j}+jy'y^iz^{j-1}).
\ees
Since  $J(f,x,h)\in A\cap J$ and
\bes
\deg_yJ(f,x,h)+\deg_zJ(f,x,h)<\deg_yf+\deg_zf,
\ees
 it follows that $J(f,x,h)=0$ or $f\in
K[x,x',y',z']\cap J$ by the choice of $f$.

Therefore, we can assume that $J(f,x,h)=0$. Then
\bes
if_{ij}z'=(j+1)f_{(i-1)(j+1)}y',\,
jf_{ij}y'=(i+1)f_{(i+1)(j-1)}z'
\ees
 and, consequently,
\bee\label{f10}
f=\sum_{i=0}^nf_i(yy'+zz')^i,
\eee
 where  $f_i\in K[x,x',y',z']$. Set $b=yy'+zz'$. Then
 $f=\sum_{i=0}^nf_ib^i$.

 Notice that $b=-\lambda-xx'-\frac{1}{4}h^2$ and we have $\{b,x\}=\{b,x'\}=\{b,h\}=0$.
We can rewrite $f$ in the form $f=\sum_{i=0}^{m}g_{i}x'^i$, where
$g_{i}\in K[x,y',z',b]$.
 Then
 \bes
 \{x,f\}=\sum_{i=0}^m g_iix'^{i-1}h=\frac{\partial f}{\partial x'}h\in J.
 \ees
Since $h\not\in J$ and   $J$ is a prime ideal it follows that
  $\frac{\partial f}{\partial x'}\in J$. This allows us to reduce the degree of $f$ in $x'$. Consequently, we can assume that  $f\in K[x,y',z',b]$.

If $f\in K[x,y',z',b]$ then we can write
\bes
f=\sum_{k=0}^n\sum_{i+j=k}f_{ij}y'^iz'^j,
\ees
 where  $f_{ij}\in K[x,b]$. Notice that
\bes
J(y',x,x')=-6y', \ J(z',x,x')=-6z', \ J(f_{ij},x,x')=0.
\ees
Consequently,
\bes
J(f,x,x')=\sum_{k=0}^n\sum_{i+j=k}f_{ij}6(-i-j)y'^iz'^j=\sum_{k=0}^n-6k\sum_{i+j=k}f_{ij}y'^iz'^j.
\ees
We also have
\bes
g=J(f,x,x')+6nf=\sum_{k=0}^{
n-1}(6n-6k)\sum_{i+j=k}f_{ij}y'^iz'^j\in J.
\ees
Using this,  we can write
\bes
0\neq f=\sum_{i+j=n}f_{ij}y'^iz'^j\in J,
\ees
 where   $f_{ij}\in
K[x,b]$. Then we can rewrite $f$ in the form $f=\sum_{k=0}^mf_{k}x^k$, where
$f_{k}=\sum_{i+j=n}g_{kij}y'^iz'^j$ and $g_{kij} \in K[b]$.  Then
\bes
\{h,f\}=\sum_{k-0}^m\{h,
f_k\}x^k+\sum_{k-0}^mf_k\{h,x^k\}=-2nf+2x\frac{\partial f}{\partial x}.
\ees
Consequently,
$2x\frac{\partial f}{\partial x}\in J$. This implies $\frac{\partial f}{\partial x}\in J$ since $x\not
\in J$. This allows us to reduce the degree of $f$ in $x$. Therefore, we can assume that $f$ does not depend on $x$ and $f_{ij}\in K[b]$.

Now we can write $f=\sum_{k=0}^mf_{k}b^k$, where
$f_{k}=\sum_{i+j=n}\alpha_{kij}y'^iz'^j$ and $\alpha_{kij} \in K$.
We have
\bes
J(f,y,h)=\sum_{k=0}^mf_kkb^{k-1}J(b,y,h)=
J(zz',y,h)\frac{\partial f}{\partial b}=-12x'z'\frac{\partial f}{\partial b}.
\ees
Consequently, $\frac{\partial f}{\partial b}\in J$ since $J(f,y,h)\in J$, $J$ is a prime ideal, and $x'z'\not\in
J$. This allows us to reduce the degree of $f$ in $b$.

Therefore, we can assume that $f$ does not depend on $b$ and $f=\sum_{i+j=n}\alpha_{ij}y'^iz'^j$, where $\alpha_{ij}\in K$, that is
$f\in K[x,x',y',z']$. Then $f\in A\cap J=0$.

This contradiction proves that $S$ is simple. $\Box$

\begin{lm}\label{l5} The center of the generic Poisson algebra
 $S=P(\mathbb M)/(c_\mathbb M-\lambda)$ is $K$.
\end{lm}
\Proof By Remark \ref{r1} we may assume that $K$ is algebraically closed. Notice that $Z(S)$ is a field since the generic Poisson algebra $S$ is simple by Lemma \ref{l4}.
Let $A=K[x,x',y,y',z,z']$. Suppose that $h$ is invertible in $S$. Then $h(f+gh)=1$ for some $f,g\in A$. Since $A\cap Ah=0$ it follows that
$f=0$ and  $gh^2=1$. Then $g\in K$ and $h^2\in K$ since $h^2\in A$. This contradiction proves that $h$ is not invertible in $S$.

Let  $a\in Z(S)$ and let $a=f+g h$,
where  $f, g\in A$. Since $J(A,b,h)\subseteq A$ for any $b\in
\{x,x',y,y'\}$ and $\{h,A\}\subseteq A$ it follows that $J(f,b,h)=0$ and
$\{h,f\}=0$.

Suppose that $f\in K[x,x',y',z']$ and $f=\sum_{i=0}^nf_ix^i$,
where  $f_i\in K[x',y',z']$. Since $J(f_i,y,h)=0$ it follows that
\bes
0=J(f,y,h)=\sum_{i=0}^n12f_iix^{i-1}z'.
\ees
This allows us to reduce the degree of $f$ in $x$. Therefore, we can assume that
$f\in
K[x',y',z']$ and
\bes
f=\sum_{l=0}^n\sum_{i+j+k=l}\alpha_{ijk}x'^iy'^jz'^k,
\ees
 where
$\alpha_{ijk}\in K$. Since $\{h,f\}=0$ it follows that
\bes
0=\{h,f\}=\sum_{l=0}^n\sum_{i+j+k=l}-2l\alpha_{ijk}x'^iy'^jz'^k.
\ees
Hence,  $f\in K$.

Recall that $gh=a-f$. This implies that $h$ is invertible or $a\in K$ since the center $Z(S)$ is a field.

Therefore, we can assume that
$f\in K[x,x',y,y',z, z']\setminus   K[x,x',y',z']$.
Suppose that $f=\sum_{i=0}^nf_iy^i$, where   $f_i\in K[x,x',y',z']$. We have $J(f_i,x,h)=0$. Then
\bes
0=J(f,x,h)=\sum_{i=0}^n-12f_iiy^{i-1}z'.
\ees
Hence $f\in
K[x,x',y',z']$. This implies that $f$ must depend on both $y$ and  $z$.

Let   $f=\sum_{ij}^nf_{ij}y^iz^j$, where  $f_{ij}\in K[x,x',y',z']$. By (\ref{f10}) we get  $f=\sum_if_i'(yy'+zz')^i$, where   $f_i'\in K[x,x',y',z']$,
since $J(f,x,h)=0$.

Let     $d=xx'+yy'+zz'$. Then  $f=\sum_if_id^i$, where  $f_i\in
K[x,x',y',z']$. Obviously, $d$ is not algebraic over $K[x,x',y',z']$. Notice that
\bes
J(x,y,h)=12z',
J(y',y,h)=J(z',y,h)=J(d,y,h)=0.
\ees
This implies that $0=J(f,y,h)=\sum_i12z'\frac{\partial
f_i}{\partial x}d^i$. Consequently, $f_i$ does not depend on $x$ and $f=\sum_if_id^i$, where  $f_i\in
K[x',y',z']$. We can rewrite $f$ in the form
$f=\sum_{l=0}^n\sum_{i+j+k=l}f_{ijk}x'^iy'^jz'^k$, where  $f_{ijk}\in
K[d]$. Using $\{h,f\}=0$ and $\{h,d\}=0$, we get
\bes
0=\{h,f\}=\sum_{l=0}^n\sum_{i+j+k=l}-2lf_{ijk}x'^iy'^jz'^k.
\ees
Consequently, $\sum_{i+j+k=l}f_{ijk}x'^iy'^jz'^k=0$ for all $l=1,2,\ldots n$.
This implies that $f=\sum_i\alpha'_{i}d^i$, where  $\alpha'_{i}\in K$.
Therefore,  $f=\sum_{i=0}^n\alpha_ih^i$, where  $\alpha_i\in K$. From here we get
\bes
a=f+gh=\alpha_0+f'h+gh,
\ees
where  $f'=\sum_i\beta_ih^i$. Since $h$ is not invertible in $S$ it follows that $a=\alpha_0\in K$.

Thus, $Z(S)=K$. $\Box$

Combining Lemmas \ref{l4} and \ref{l5}, we get the following theorem.

\begin{theor} \label{t2} Let $\mathbb M={\mathbb M}(\alpha,\beta,\gamma)$ be a seven dimensional simple Malcev algebra over a field $K$ of characteristic zero with nonzero parameters
$\alpha, \beta,\gamma\in K$.  Let $c_{\mathbb M}=c_{\mathbb M}(\alpha,\beta,\gamma)$ be the standard Casimir element of $P(\mathbb M)$ and let $0\neq \lambda\in K$.  Then the generic Poisson algebra
 $S=P(\mathbb M)/(c_\mathbb M-\lambda)$ is central and simple.
\end{theor}

\end{document}